\documentclass[final]{siamltex}

\usepackage{graphicx}
\usepackage{amssymb}
\usepackage{epstopdf}
\usepackage{cite}
\usepackage{rotfig_tikz,pgfplots}
\usepackage[caption=false]{subfig}

\pgfplotsset{width=10cm,compat=1.9}

\usepackage{booktabs}

\newcommand{\cplxs}{\mbox{$\mathbb{C}$}}

\newcommand{\cnxn}{\mbox{$\cplxs^{n \times n}$}}
\newcommand{\absval}[1]{\mbox{$\mid\!#1\!\mid$}}
\newcommand{\norm}[1]{\mbox{$\parallel\!#1\!\parallel$}}

\newcommand{\spn}[1]{\mbox{{\rm span}$\left\{#1\right\}$}}
\newcommand{\eee}{\mbox{$\mathcal{E}$}}

\newcommand{\cue}{\mbox{$\mathcal{Q}$}}

\title{Pole-swapping algorithms for alternating and palindromic eigenvalue 
problems\thanks{This research was partially supported by 
the Research Council KU Leuven, project 
C14/16/056 (Inverse-free Rational Krylov Methods: Theory and Applications).}}
\author{Thomas Mach\footnotemark[3]\and 
Thijs Steel\footnotemark[6]
Raf Vandebril\footnotemark[4]\and
David~S.~Watkins\footnotemark[5]}

\begin{document}

\date{today}
\maketitle

\renewcommand{\thefootnote}{\fnsymbol{footnote}}

\footnotetext[3]{Department of Mathematical Sciences, Kent State University, Ohio
  (\texttt{tmach1@kent.edu}).}%
\footnotetext[6]{Department of Computer Science, KU Leuven, Belgium 
(\texttt{thijs.steel@cs.kuleuven.be}).}
\footnotetext[4]{Department of Computer Science, KU Leuven,
  Belgium (\texttt{raf.vandebril@cs.kuleuven.be}).}
\footnotetext[5]{Department of Mathematics, Washington State University
   (\texttt{watkins@math.wsu.edu}).}

\renewcommand{\thefootnote}{\arabic{footnote}}

\begin{abstract}
Pole-swapping algorithms are generalizations of bulge-chasing algorithms
for the generalized eigenvalue problem.   Structure-preserving pole-swapping algorithms for the  
palindromic and alternating eigenvalue problems, which arise in control theory, are derived.   
A refinement step that guarantees backward stability of the algorithms is included.  This 
refinement can also be applied to bulge-chasing algorithms that had been introduced 
previously, thereby guaranteeing their backward stability in all cases.  
\end{abstract}



\section{Introduction}
Francis's implicitly-shifted QR algorithm \cite{Fra61b,Wat11} is still the 
standard tool for computing the eigenvalues of a small to medium-sized  non-Hermitian matrix
$A\in\cnxn$.  Eigenvalue problems often arise naturally as generalized eigenvalue problems for
a pencil $A - \lambda B$, and for these the Moler-Stewart variant \cite{MolSte73} 
of Francis's algorithm, commonly called the QZ algorithm,  is used.   The Francis and Moler-Stewart 
algorithms are prime examples of what we call \emph{bulge-chasing} algorithms.  

In recent years some generalizations of the QZ algorithm have been proposed, e.g.\  
\cite{VanWat12g} and, more generally, the rational QZ (RQZ) algorithm of Camps, Meerbergen, and Vandebril \cite{CaMeVa19a}, which
is the first example of what we call a \emph{pole-swapping} algorithm.  
This work has been extended in various directions in \cite{Cam19,CaMaVaWa19,CaMeVa19b}.
 
In this paper we extend in another direction, 
introducing structure-preserving pole-swapping algorithms for two classes of 
structured generalized eigenvalue problems that arise in optimal control theory \cite{Meh91}, 
namely palindromic and alternating eigenvalue problems. 
Kressner, Schr\"{o}der, and Watkins \cite{KrScWa09} proposed 
structure-preserving bulge-chasing (QZ-like) algorithms for palindromic and alternating eigenvalue problems.
We show that our structured pole-swapping algorithms are generalizations of these bulge-chasing algorithms.  
Our algorithms include a refinement step, which can also be incorporated into the algorithm of \cite{KrScWa09},
to ensure backward stability.

We have opted to cover only the complex case.  If the matrices $A$ and $B$ happen to be real, the algorithms 
introduced here will not preserve the real structure.   It is possible, and not difficult, to build pole-swapping algorithms
that preserve real structure, as was shown in \cite{Cam19,CaMeVa19b}.  We could have incorporporated
that extension in this paper at the expense of making it longer and less readable.  

\section{Basic definitions and facts}

In this paper we will refer sometimes to a pencil $A - \lambda B$ and other times to a pair $(A,B)$.  
Either way, we are talking about the same object.   The pencil $A - \lambda B$ is called \emph{singular} if 
the polynomial $\det(A - \lambda B)$ is identically zero.  Otherwise the pencil is \emph{regular}.   Throughout the paper 
we make the blanket assumption that we are dealing with a regular pencil.  

The pair $(A,B)$ is called \emph{palindromic} (or \emph{$*$-palindromic}) if $A$ and $B$ are related by 
\begin{equation}\label{eq:palindromic}
A^{*} = B.
\end{equation}
The pair $(A,B)$ is \emph{alternating}  (or \emph{even}) if 
\begin{equation}\label{eq:alternating}
A^{*} = A \quad\mbox{and}\quad B^{*} = - B.
\end{equation}
These two structures are equivalent in principle, since the generalized Cayley transform 
$(A,B) \to (A+B,A-B)$ maps a palindromic pair to an alternating pair and vice versa. 

Each of these structures exhibits a spectral symmetry.  In the palindromic case, 
$\lambda$ is an eigenvalue if and only if $1/\overline{\lambda}$ is.  This is shown by writing down the equation
$Ax = \lambda Bx$, taking the conjugate transpose, and applying (\ref{eq:palindromic}).   Clearly $\lambda = 
1/\overline{\lambda}$ if and only if $\lambda$ is on the unit circle.  Eigenvalues on the unit circle need not occur in 
pairs, but those off the unit circle must appear in $(\lambda,1/\overline{\lambda})$ pairs, one inside and one outside the 
unit circle.

In the alternating case, $\lambda$ is an eigenvalue if and only if $-\overline{\lambda}$ is, as can be shown in the same
way as for the palindromic case.   Since $\lambda = -\overline{\lambda}$ if and only if $\lambda$ is on the imaginary
axis,  eigenvalues on the imaginary axis need not occur in pairs, but those off of the imaginary axis must occur in 
$(\lambda,-\overline{\lambda})$ pairs, one on each side of the imaginary axis.  

The spectral symmetry in alternating pencils is exactly the same as that possessed by Hamiltonian matrices.  
Recall that a matrix $H\in\cplxs^{2m\times 2m}$ is called \emph{Hamiltonian} if $JH$ is Hermitian, where
\begin{displaymath}
J = \left[\begin{array}{cc} & I_{m} \\ -I_{m} & \end{array}\right].
\end{displaymath} 
It is well known \cite{Meh91} that the continuous-time
linear-quadratic control problem can be solved by computing the eigensystem
of a Hamiltonian matrix.  One can equally well formulate the problem as an eigenvalue problem for an alternating 
pencil.  In fact, the Hamiltonian eigenvalue problem $H - \lambda I$ is clearly equivalent to the alternating 
eigenvalue problem $JH - \lambda J$.    Hamiltonian eigenvalue problems are also sometimes studied in the 
guise of Hamiltonian/skew-Hamiltonian pencils, but note that the Hamiltonian/skew-Hamiltonian pencil 
$H - \lambda S$ is equivalent to the alternating pencil $JH - \lambda JS$.   Alternating pencils are discussed
in various guises and contexts in 
\cite{MehWat01,MehWat02,ApMeWa02,MaMaMeMe06a,MaMaMeMe06b}, for example.

The spectral symmetry in palindromic pencils is the same as that possessed by symplectic matrices.  
The discrete-time linear-quadratic optimal control problem can be solved by computing the eigensystem
of a symplectic matrix or pencil \cite{Meh91}, 
which can also be formulated as a palindromic eigenvalue problem
\cite{MaMaMeMe06b,MaMaMeMe09,KrScWa09}.

In \cite{KrScWa09} we considered pairs $(A,B)$ for which both $A$ and $B$ have anti-Hessenberg form:
\begin{displaymath}
\parbox{3.2cm}{
\begin{tikzpicture}[scale=1.66,y=-1cm]
\draw (-.15,-.1) -- (-.2,-.1) -- (-.2,1.5) -- (-.15,1.5);
\draw (1.55,-.1) -- (1.6,-.1) -- (1.6,1.5) -- (1.55,1.5);
\foreach \j in {0,...,7}{
   \foreach \i in {\j,...,7}{\node at (\i/5,1.4-\j/5)
     [align=center,scale=1.0]{$\times$};}}
\foreach \j in {0,...,6}{\node at (\j/5,1.4-\j/5-.2)
     [align=center,scale=1.0]{$\times$};}
\end{tikzpicture}
},
\end{displaymath}
shown here in the case $n=8$.  We will consider the same form here.  We could study instead 
the equivalent pair $(FA,FB)$, where $F$ is the \emph{flip} or 
\emph{anti-identity} matrix.  $FA$ and $FB$ are both 
upper Hessenberg, but the anti-Hessenberg form is more convenient for the special structures that 
we are considering here.  

This form is admittedly quite special.  Given a general pair $(A,B)$ that is either palindromic or alternating,
there is no known efficient method for transforming the pair to anti-Hessenberg form while preserving the 
relevant structure, except in special cases.   
For example, if the pencil corresponds to a control system that has either a single
input or a single output, such a reduction is possible \cite{KrScWa09}.  

Following Camps, Meerbergen, and Vandebril \cite{CaMeVa19a} we associate $n-1$ \emph{poles} with the 
anti-Hessenberg pair $(A,B)$.  For $k=1$, \ldots, $n-1$ the \emph{pole in position $(n-k,k)$} is the 
ratio $\sigma_{k} = a_{n-k,k}/b_{n-k,k}$.    We assume that for each $k$, $\absval{a_{n-k,k}} + \absval{b_{n-k,k}} > 0$, 
since otherwise it would be possible to reduce the problem immediately to two or more smaller eigenvalue problems.  
Thus every $\sigma_{k}$ is well defined (but might equal~$\infty$).  

In the two structured cases that we are considering, the poles exhibit the same symmetry as the eigenvalues do.
In the palindromic case we have 
$a_{n-k,k} = \overline{b}_{k,n-k}$, and it follows that 
\begin{displaymath}
\sigma_{k} = 1/\overline{\sigma}_{n-k},\quad k=1,\, \ldots,\, n-1.   
\end{displaymath}
In the case of even $n$, there is one unpaired pole $\sigma_{n/2}$, which must satisfy 
$\absval{\sigma_{n/2}} = 1$.  
In the alternating case we have 
\begin{displaymath}
\sigma_{k} = -\overline{\sigma}_{n-k}, \quad k=1,\, \ldots,\, n-1.
\end{displaymath}   
If $n$ is even, there is 
one unpaired pole $\sigma_{n/2}$, which must lie on the imaginary axis.  

\section{Operations on anti-Hessenberg pairs}\label{sec:operations}

Introducing terminology that we have used in some of our recent work 
\cite{AuMaRoVaWa18,AuMaRoVaWa18g,AuMaRoVaWa19,AuMaVaWa15},  we define a \emph{core transformation}  
(or \emph{core} for short) to be a unitary matrix that acts only on two adjacent rows/columns, for example,
\begin{displaymath}
Q_{3} = \left[\begin{array}{ccccc}
1 & & & & \\ & 1 & & & \\ & & {*} & {*} & \\  & & {*} & {*} & \\  & & & & 1
\end{array}\right],
\end{displaymath}
where the four asterisks form a $2 \times 2$ unitary matrix.  Givens rotations are examples of core transformations.
Our core transformations always have subscripts that tell where the action is:  $Q_{j}$ acts on rows/columns
$j$ and $j+1$.

In \cite{CaMeVa19a} two types of operations on upper Hessenberg pencils were introduced.  These are 
unitary equivalence transformations, which we called \emph{moves} of types I and II in \cite{CaMaVaWa19}. 
Obviously we can do the same sorts of moves on anti-Hessenberg  pencils, but if we want to preserve the 
palindromic or alternating structure, we need to use special unitary equivalence transformations, namely 
congruences  $A - \lambda B \to Q^{*}(A - \lambda B)Q$.  

\subsection*{Move of Type I}  
This move replaces the pole $\sigma_{1}$  (located at position $(n-1,1)$) by any other value $\rho$.  
At the same time, since the symmetry must be preserved, the pole $\sigma_{n-1}$ is changed appropriately.
In the palindromic case, $\sigma_{n-1} = 1/\overline{\sigma}_{1}$, and it gets changed to $1/\overline{\rho}$.
In the alternating case, $\sigma_{n-1} = -\overline{\sigma}_{1}$, and it gets changed to $-\overline{\rho}$.

To see how this is done, we at first ignore the need to preserve structure and consider how we would 
insert a pole $\rho$ at position $(n-1,1)$.  Because of the anti-Hessenberg form,  the vector $(A - \rho B)e_{1}$ 
consists of zeros, except for the last two entries.\footnote{Here and in what follows, the notation $A - \rho B$ is shorthand
for $\beta A - \alpha B$, where $\alpha$ and $\beta$ are any scalars satisfying $\rho = \alpha/\beta$.  This allows us to
include the case $\rho = \infty$ by taking $\beta = 0$.}   Therefore a core transformation $Q_{n-1}$, acting on rows
$n-1$ and $n$, can be constructed so that $Q_{n-1}^{*}$ zeros out the entry in position $n-1$, that is, 
\begin{displaymath}
Q_{n-1}^{*}(A - \rho B)e_{1} = \gamma\,e_{n} 
\end{displaymath}
for some nonzero $\gamma$.   Now let $\tilde{A}  - \lambda \tilde{B} 
= Q_{n-1}^{*}(A - \lambda B)$.  This new pencil has the pole $\rho$ in position $(n-1,1)$, as desired, since 
$\tilde{a}_{n-1,1} - \rho\tilde{b}_{n-1,1} = 0$.    This is exactly the move of type I described in \cite{CaMaVaWa19} and
earlier in \cite{CaMeVa19a}, turned over for the anti-Hessenberg case.

But this transformation does not preserve the symmetry of the pencil.  What is needed is a congruence:  
$\hat{A} - \lambda \hat{B} = Q_{n-1}^{*}(A - \lambda B)Q_{n-1}$.  The right multiplication by $Q_{n-1}$ does not
affect what happens in the lower left-hand corner of the pencil, so the pole in position $(n-1,1)$ is $\rho$, as shown
above.  But the right multiplication by $Q_{n-1}$ does affect the pole $\sigma_{n-1}$ at position $(1,n-1)$, changing it
(by symmetry) to $1/\overline{\rho}$ in the palindromic case and $-\overline{\rho}$ in the alternating case.

\subsection*{Move of Type II}
This move swaps two adjacent poles.  If we delete the $n$th row and column from the anti-Hessenberg pencil
$A - \lambda B$, we get an anti-triangular pencil $A_{\pi} -\lambda B_{\pi}$, which we call the \emph{pole pencil} 
because its eigenvalues are exactly the poles of $A - \lambda B$.  Swapping two adjacent poles in $A - \lambda B$
is the same as swapping two eigenvalues in $A_{\pi} - \lambda B_{\pi}$, and there are well-known procedures for 
doing this \cite{BaiDem93}, \cite{CaMaVaWa19}, \cite{KagPor96a,KagPor96b}, 
\cite{VanD81}, \cite[\S\S~4.8,\,6.6]{Wat07}.   

Suppose 
we want to swap two adjacent poles $\sigma_{k-1}$ and $\sigma_{k}$
located at positions $(n-k+1,k-1)$ and $(n-k,k)$.  Temporarily ignoring symmetry, this can
be accomplished by a transformation $Q_{n-k}^{*}(A - \lambda B)Z_{k-1}$.   All of the action is
in the subpencil 
\begin{equation}\label{eq:2pencil}
\left[\begin{array}{cc} 0 & a_{n-k,k} \\ a_{n-k+1,k-1} & a_{n-k+1,k}\end{array}\right]
 - \lambda \left[\begin{array}{cc} 0 & b_{n-k,k} \\ b_{n-k+1,k-1} & b_{n-k+1,k}\end{array}\right].
\end{equation}
This is the principal subpencil of the bulge pencil that contains the two poles $\sigma_{k-1} = 
a_{n-k+1,k-1}/b_{n-k+1,k-1}$ and $\sigma_{k} = a_{n-k,k}/b_{n-k,k}$.

In order to preserve the symmetry we will also have 
to swap the poles $\sigma_{n-k}$ and $\sigma_{n-k+1}$  
in positions $(k-1,n-k+1)$ and $(k,n-k)$.   The appropriate transformation is, by symmetry,
$Z_{k-1}^{*}(A - \lambda B)Q_{n-k}$.  This can be done provided the two transformations do not
interfere with each other, i.e. $Q_{n-k}$ and $Z_{k-1}$ act on non-overlapping columns.  This is the case
if $k < n-k$ or $n-k+1 < k-1$, that is, $k < n/2$ or $k > n/2 + 1$.
The total transformation is $Q^{*}(A - \lambda B)Q$, where $Q = Z_{k-1}Q_{n-k} = Q_{n-k}Z_{k-1}$.

\subsection*{A hint at where we are heading}

In \cite{CaMeVa19a,CaMaVaWa19}  it was shown how moves of types I and II can be used to build algorithms 
for computing eigenvalues of an upper Hessenberg pencil.  In the simplest case a 
shift $\rho$  is chosen, and a move of 
type I is used to insert it as a pole at the top of the pencil.  
Then a sequence of moves of type II is used to exchange the pole $\rho$ downward until it gets to the bottom
of the pencil.  Then it is removed from the bottom 
(replaced by some new pole $\sigma_{n}$) by a move of type I.   This procedure
can be shown \cite{CaMeVa19a} to be a generalization of (one iteration of) the QZ algorithm on a Hessenberg-triangular 
pencil.

In our current scenario the matrices are anti-Hessenberg, not Hessenberg, but that is a trivial difference.  More
importantly, in the moves described here, everything is doubled up for preservation of 
structure.  If we introduce a pole $\rho$ at one end of the pencil, we must simultaneously introduce a pole 
$\tilde{\rho}$ ($=1/\overline{\rho}$ or $-\overline{\rho}$, depending on the structure)  at the other end.  Now if we want to move $\rho$ from one end of the pencil to the other by moves of type II, we must simultaneously move $\tilde{\rho}$ in the opposite 
direction.  There comes a point in the middle where we have to swap $\rho$ and $\tilde{\rho}$, so that they can continue
their journey to the opposite end of the pencil.  This requires
a special symmetric version of the move of type II.

\subsection*{Move of Type IIo}  

Suppose the dimension $n$ is odd.  Then there is an even number of poles $\sigma_{1}$, \ldots, $\sigma_{n-1}$. 
The two poles in the middle are $\sigma_{k-1}$ and $\sigma_{k}$, where $k = (n+1)/2$.  
These are the two eigenvalues of (\ref{eq:2pencil}) with $k=(n+1)/2$.  
In the interest of simplicity and non-proliferation of notation, we rename this subpencil 
\begin{displaymath}
A - \lambda B = 
\left[\begin{array}{cc} & \alpha_{1} \\ \alpha_{2} & \alpha_{21} \end{array}\right]
- \lambda 
\left[\begin{array}{cc} & \beta_{1} \\ \beta_{2} & \beta_{21} \end{array}\right].
\end{displaymath}
We are temporarily re-assigning the symbols $A$ and $B$ to stand for the submatrices that are our current focus.  
This little pole pencil has the same structure as the original pencil, either palindromic or alternating, but we will 
ignore the structure at first.   The eigenvalues are $\alpha_{2}/\beta_{2}$ and $\alpha_{1}/\beta_{1}$, and we would like to
swap them.  If we flip the rows and columns we get 
\begin{displaymath}
FAF - \lambda FBF = 
\left[\begin{array}{cc} \alpha_{21} & \alpha_{2} \\ \alpha_{1} &  \end{array}\right]
- \lambda 
\left[\begin{array}{cc} \beta_{21} & \beta_{2} \\ \beta_{1} &  \end{array}\right],
\end{displaymath}
which has the poles in the desired positions but the ``wrong'' triangularity.  We have to fix this.  We
find it convenient to work with the partially flipped form 
\begin{displaymath}
AF - \lambda BF = 
\left[\begin{array}{cc} \alpha_{1} &   \\ \alpha_{21} & \alpha_{2} \end{array}\right]
- \lambda 
\left[\begin{array}{cc} \beta_{1} &  \\ \beta_{21} & \beta_{2} \end{array}\right],
\end{displaymath}
which is easier to study because the matrices are triangular, not anti-triangular.  

Next we set up and solve a Sylvester equation to diagonalize the pencil.  Specifically, we will find 
unit lower triangular matrices 
\begin{displaymath}
X = \left[\begin{array}{cc} 1 &  \\ x & 1  \end{array}\right]
\quad\mbox{and}\quad
Y = \left[\begin{array}{cc} 1 &  \\ y & 1 \end{array}\right]
\end{displaymath}
such that 
\begin{equation}\label{eq:sylv}
(AF)X = Y(\check{A}F) \quad\mbox{and}\quad (BF)X = Y(\check{B}F),
\end{equation}
where $\check{A}$ and $\check{B}$ are anti-diagonal matrices with the same anti-diagonals as $A$ and $B$, 
respectively.  Writing the first of these equations out in detail, we have 
\begin{displaymath}
\left[\begin{array}{ccc} \alpha_{1} &  \\
\alpha_{21} & \alpha_{2} \end{array}\right]
\left[\begin{array}{ccc} 1 &  \\ x & 1  \end{array}\right] =
\left[\begin{array}{ccc} 1 & \\ y & 1 \end{array}\right]
\left[\begin{array}{ccc} \alpha_{1} &     \\
& \alpha_{2} \end{array}\right],
\end{displaymath}
and similarly for the $B$ equation.  This is a system of two linear equations in two unknowns
\begin{displaymath}
\alpha_{21} + \alpha_{2}x  = y \alpha_{1}, \quad \beta_{21} + \beta_{2}x = y \beta_{1}
\end{displaymath}
or 
\begin{equation}\label{eq:sylv0}
\left[\begin{array}{cc} 
\alpha_{1} & \alpha_{2}  \\ \beta_{1} & \beta_{2}
\end{array}\right]
\left[\begin{array}{c} 
\phantom{-}y \\ -x
\end{array}\right]  =
\left[\begin{array}{c} \alpha_{21} \\ \beta_{21}  \end{array}\right].
\end{equation}
This has a unique solution if and only if $\alpha_{1}\beta_{2} - \beta_{1}\alpha_{2} \neq 0$, i.e.\ the
poles are distinct.  

Thus, assuming distinct poles, (\ref{eq:sylv}) has a unique solution, which we can easily and stably compute.  
Rewriting (\ref{eq:sylv}) we have the equivalence 
\begin{equation}\label{eq:nonuneq2}
 FY^{-1}(A - \lambda B)(FX) = F(\check{A} - \lambda\check{B})F = 
\left[\begin{array}{cc} & \alpha_{2}  \\  \alpha_{1} & \end{array}\right]
- \lambda
\left[\begin{array}{ccc} & \beta_{2}   \\ \beta_{1} & \end{array}\right],
\end{equation} 
which (still) has the poles in the desired locations,  but this equivalence is not unitary.  To make 
a unitary equivalence we will introduce some QR decompositions.  Since we are going to work with 
triangular matrices, we again remove the $F$ from the left-hand side to obtain 
\begin{displaymath}
Y^{-1}(A - \lambda B)(FX) = (\check{A} - \lambda\check{B})F = 
\left[\begin{array}{cc} \alpha_{1} & \\ & \alpha_{2}  \end{array}\right]
- \lambda
\left[\begin{array}{ccc} \beta_{1} &  \\ & \beta_{2} \end{array}\right].
\end{displaymath} 
 Let 
 \begin{equation}\label{eq:2qrd}
 FX = QR \quad\mbox{and}\quad Y = PS,
 \end{equation}
 where $Q$ and $P$ are unitary, and $R$ and $S$
 are upper triangular with positive entries on the main diagonal.   Then 
 \begin{displaymath}
 P^{*}(A - \lambda B)Q = S(\check{A}F -\lambda \check{B}F)R^{-1} = T_{A} -\lambda T_{B},
 \end{displaymath}
 where $T_{A}$ and $T_{B}$ are upper triangular:
 \begin{displaymath}
 T_{A} = S(\check{A}F) R^{-1} = \left[\begin{array}{cc} 
 s_{11}\alpha_{1} r_{11}^{-1} &  {*} \\ & s_{22}\alpha_{2}r_{22}^{-1} 
 \end{array}\right],
 \end{displaymath}
 \begin{displaymath}
 T_{B} = S(\check{B}F) R^{-1} = \left[\begin{array}{cc} 
 s_{11}\beta_{1} r_{11}^{-1} &  {*} \\ & s_{22}\beta_{2}r_{22}^{-1} 
 \end{array}\right].
 \end{displaymath}
 Since we want anti-triangular matrices, we now restore the $F$ on the left to obtain
 \begin{equation}\label{eq:finaleq}
 FP^{*}(A - \lambda B) Q =  
\left[\begin{array}{cc} 
 & s_{22}\alpha_{2} r_{22}^{-1} \\
 s_{11}\alpha_{1}r_{11}^{-1} & {*}
 \end{array}\right] - \lambda \left[\begin{array}{cc} 
  & s_{22}\beta_{2} r_{22}^{-1} \\
 s_{11}\beta_{1}r_{11}^{-1} & {*}
 \end{array}\right].  
 \end{equation}
 This unitary equivalence gives the right anti-triangular form with the poles 
 $\alpha_{1}/\beta_{1}$ and $\alpha_{2}/\beta_{2}$ in the desired positions.  

So far we have assumed no special relationship between $A$ and $B$.   To finish the story we must 
show that in our two structured cases  (\ref{eq:finaleq}) is a congruence, that is $PF = Q$.  
 In preparation
for this we go back to (\ref{eq:sylv}) and take complex conjugates of both equations to obtain 
\begin{displaymath}
X^{*}FA^{*} = F\check{A}^{*}Y^{*}, \quad X^{*}FB^{*} = F\check{B}^{*}Y^{*}.  
\end{displaymath}
Further simple manipulations yield 
\begin{equation}\label{eq:sylvstar}
(A^{*}F)(FY^{-*}F) = (FX^{-*}F)(\check{A}^{*}F), \quad (B^{*}F)(FY^{-*}F) = (FX^{-*}F)(\check{B}^{*}F).
\end{equation}
Now consider what these equations look like in the alternating case $A^{*} = A$, $B^{*} = - B$ 
(which also implies $\check{A}^{*} = \check{A}$ and $\check{B}^{*} = - \check{B}$).  Clearly
we have
\begin{equation}\label{eq:sylvstarstruc}
(AF)(FY^{-*}F) = (FX^{-*}F)(\check{A}F), \quad (BF)(FY^{-*}F) = (FX^{-*}F)(\check{B}F).
\end{equation}
Now consider the palindromic case $B = A^{*}$.  Inserting this equation into (\ref{eq:sylvstar}), we again 
get (\ref{eq:sylvstarstruc}), just as in the alternating case.  Thus, in either case, (\ref{eq:sylvstarstruc}) holds.  
Noting that $FY^{-*}F$ and $FX^{-*}F$ are both 
unit lower triangular, and comparing (\ref{eq:sylvstarstruc}) with (\ref{eq:sylv}), we see that the pair 
$(FY^{-*}F,FX^{-*}F)$ is a solution of (\ref{eq:sylv}).  By uniqueness of the solution we deduce that 
$(X,Y) = (FY^{-*}F,FX^{-*}F)$, and in particular $X = FY^{-*}F$.  Writing this
as $FX = Y^{-*}F$ and inserting the decompositions $FX = QR$ and $Y= PS$, we obtain
\begin{displaymath}
QR = PS^{-*}F = (PF)(FS^{-*}F). 
\end{displaymath}  
By uniqueness of the QR decomposition, 
\begin{displaymath}
Q = PF.
\end{displaymath}
Making this substitution into (\ref{eq:finaleq}), we obtain
\begin{displaymath}
Q^{*}(A - \lambda B)Q = FT_{A} - \lambda FT_{B} = \hat{A} - \lambda\hat{B},
\end{displaymath}
where

\begin{displaymath}
\hat{A} = \left[\begin{array}{cccc}  & \alpha_{2}s_{22}r_{22}^{-1} \\
\alpha_{1}s_{11}r_{11}^{-1} & {*} 
\end{array}\right]
\quad\mbox{and}\quad 
\hat{B} = \left[\begin{array}{cccc}  & \beta_{2}s_{22}r_{22}^{-1} \\
\beta_{1}s_{11}r_{11}^{-1} & {*} 
\end{array}\right].
\end{displaymath}

The actual computation is quite simple.  From (\ref{eq:2qrd}) we see that we need to build  $Q$ 
satisfying $FX = QR$.  This means we need
\begin{equation}\label{eq:qfirstcol}
\left[\begin{array}{c} x \\ 1\end{array}\right] = Q\left[\begin{array}{c} r_{11} \\ 0 \end{array}\right], 
\end{equation}
where $x$ is obtained by solving (\ref{eq:sylv0}).   In other words, we need a unitary $Q$ whose 
first column is proportional to $\left[\begin{array}{cc} x & 1 \end{array}\right]^{T}$.
Partially solving (\ref{eq:sylv0}) by Cramer's rule, we obtain 
\begin{displaymath}
x = \frac{\beta_{1}\alpha_{12} - \alpha_{1}\beta_{12}}{\alpha_{1}\beta_{2} - \beta_{1}\alpha_{2}}.
\end{displaymath}
We don't actually need $\left[\begin{array}{cc} x & 1 \end{array}\right]^{T}$ in (\ref{eq:qfirstcol}), as
any vector proportional to $\left[\begin{array}{cc} x & 1 \end{array}\right]^{T}$ will suffice, so 
we use instead 
\begin{displaymath}
v = \left[\begin{array}{c} v_{1} \\ v_{2}\end{array}\right] = 
\left[\begin{array}{c} \beta_{1}\alpha_{12} - \alpha_{1}\beta_{12}  \\ 
\alpha_{1}\beta_{2} - \beta_{1}\alpha_{2}
\end{array}\right].
\end{displaymath}
Thus $Q$ is obtained by 
\begin{displaymath}
r = \norm{v}_{2}, \qquad c = \frac{v_{1}}{r}, \qquad 
s = \frac{v_{2}}{r},
\end{displaymath}
and 
\begin{displaymath}
Q = \left[\begin{array}{cc}  c & -\overline{s} \\ 
s & \phantom{-}\overline{c}\end{array}\right].
\end{displaymath}
We remark that in both palindromic and alternating cases, $s$ is a real number.  

In terms of the original 
(big) matrices $A$ and $B$, the congruence is 
$$Q_{k-1}^{*}(A - \lambda B)Q_{k-1},$$ 
where $Q_{k-1}$ is 
a core transformation built from the $2\times 2$ unitary matrix $Q$.

We emphasize that this procedure succeeds if the two poles are distinct.  Of course there is nothing
to be gained by interchanging two poles that are equal.  In our application below we will use a move of this 
type to interchange two shifts $\rho$ and $\tilde{\rho}$.  In the palindromic (resp.\ alternating) case, 
$\tilde{\rho} = 1/\overline{\rho}$  (resp.\ $\tilde{\rho} = - \overline{\rho}$),
and the equation $\rho \neq \tilde{\rho}$ just means that $\rho$ does not lie on the unit circle 
(resp.\ imaginary axis). 

\subsection*{Refinement of a move of type IIo}

After the move of type IIo, the resulting $2 \times 2$ pole pencil has the form 
\begin{displaymath}
\hat{A} - \lambda\hat{B} =  
\left[\begin{array}{cc} 0 & \hat{\alpha}_{2} \\ \hat{\alpha}_{1} & \hat{\alpha}_{12} \end{array}\right]
- \lambda 
\left[\begin{array}{cc} 0 & \hat{\beta}_{2} \\ \hat{\beta}_{1} & \hat{\beta}_{12} \end{array}\right]
\end{displaymath}
in principle.  In practice the numbers that are supposed to be zero will not be exactly zero because of
roundoff errors in the computation.  
If those numbers are small enough (a modest multiple of the unit roundoff), they can 
simply be set to zero without compromising stability.  If, on the other hand, they are not small enough,
a correction step can make them smaller.  Now let's rewrite the pencil to take the nonzero entries into account.  
At the same time we will recycle the notation, leaving off the hats for simplicity.  We have
\begin{displaymath}
A - \lambda B  = \left[\begin{array}{cc} \epsilon & \alpha_{2} \\ \alpha_{1} & \alpha_{12} \end{array}\right]
- \lambda \left[\begin{array}{cc} \eta & \beta_{2} \\ \beta_{1} & \beta_{12} \end{array}\right],
\end{displaymath}
where $\absval{\epsilon}$ and $\absval{\eta}$ are tiny but not small enough to be set to zero.  

The correction step is explained in detail, and in greater generality, in Section~\ref{sec:refine}.  Here we provide
a brief description.  We look for tiny corrections $x$ and $y$ such that 
\begin{displaymath}
\left[\begin{array}{cc} 1 & y \\ & 1 \end{array}\right] 
\left[\begin{array}{cc} \epsilon & \alpha_{2} \\ \alpha_{1} & \alpha_{12} \end{array}\right]
\left[\begin{array}{cc} 1 &  \\  x & 1 \end{array}\right] = 
\left[\begin{array}{cc} 0 & \check{\alpha}_{2} \\ \check{\alpha}_{1} & \check{\alpha}_{12} \end{array}\right], 
\end{displaymath}
and similarly for $B$.  This yields a pair of equations 
\begin{displaymath}
\epsilon + y\alpha_{1} + \alpha_{2}x + y\alpha_{12}x = 0, \quad 
\eta + y\beta_{1} + \beta_{2}x + y\beta_{12}x = 0, 
\end{displaymath}
which can be simplified by deleting the insignificant quadratic terms to yield the linear equations
\begin{equation}\label{eq:sylv1}
\left[\begin{array}{cc} \alpha_{1} & \alpha_{2} \\ \beta_{1} & \beta_{2} \end{array}\right]
\left[\begin{array}{c} y \\ x\end{array}\right] = - \left[\begin{array}{c} \epsilon \\ \eta \end{array}\right].
\end{equation}
Since the poles $\alpha_{1}/\beta_{1}$ and $\alpha_{2}/\beta_{2}$ are distinct,  this system has a unique solution.
It is easy to show that in both the alternating and the palindromic cases, $y = \overline{x}$, so the contemplated 
transformation is in fact a congruence.  To make a unitary congruence we do a QR decomposition
\begin{displaymath}
\left[\begin{array}{cc} 1 & \\ x & 1 \end{array}\right] = QR.
\end{displaymath}
Because $x$ is tiny, $Q$ is close to the identity matrix.  The orthogonal congruence by $Q$ is the desired correction:
\begin{displaymath}
Q^{*} \left( \left[\begin{array}{cc} \epsilon & \alpha_{2} \\ \alpha_{1} & \alpha_{12} \end{array}\right]
- \lambda \left[\begin{array}{cc} \eta & \beta_{2} \\ \beta_{1} & \beta_{12} \end{array}\right] \right) Q =
\left[\begin{array}{cc} \hat{\epsilon} & \hat{\alpha}_{2} \\ \hat{\alpha}_{1} & \hat{\alpha}_{12} \end{array}\right]
- \lambda \left[\begin{array}{cc} \hat{\eta} & \hat{\beta}_{2} \\ \hat{\beta}_{1} & \hat{\beta}_{12} 
\end{array}\right],
\end{displaymath}
where $\hat{\epsilon}$ and $\hat{\eta}$ are now normally small enough to be set to zero.  
In the unlikely event that they are not, the correction step can be repeated.

Again the computation is quite simple.  Solving for $x$ in (\ref{eq:sylv1}),  we obtain 
\begin{displaymath}
x = \frac{\beta_{1}\epsilon - \alpha_{1}\eta}{\alpha_{1}\beta_{2} - \beta_{1}\alpha_{2}}.
\end{displaymath}
We need a unitary $Q$ whose first column is proportional to $\left[\begin{array}{cc} 1 & x\end{array}\right]^{T}$
or equivalently
\begin{displaymath}
v = \left[\begin{array}{c} v_{1} \\ v_{2}\end{array}\right] = 
\left[\begin{array}{c}  \alpha_{1}\beta_{2} - \beta_{1} \alpha_{2} \\ \beta_{1}\epsilon - \alpha_{1}\eta \end{array}\right]. 
\end{displaymath}
Thus we can take 
\begin{displaymath}
Q = \left[\begin{array}{cc} c & -\overline{s} \\ s & \phantom{-}\overline{c}\end{array}\right],
\end{displaymath}
where $r = \norm{v}_{2}$, $c = v_{1}/r$, and $s = v_{2}/r$.

\subsection*{Move of Type IIe}

Now consider the case when $n$ is even.  There is an odd number of poles $\sigma_{1}$, \ldots, $\sigma_{n-1}$.  
The three poles in the middle are $\sigma_{k-1}$, $\sigma_{k}$, and $\sigma_{k+1}$, where $k = n/2$.  Recall that 
the pole $\sigma_{k}$ is a special unpaired pole, while $\sigma_{k-1}$ and $\sigma_{k+1}$ are linked by the
structure (palindromic or alternating).  We need a move that swaps $\sigma_{k-1}$ with $\sigma_{k+1}$, while
leaving $\sigma_{k}$ where it is.  

All of the action takes place in a $3 \times 3$ subpencil 
\begin{displaymath}
A - \lambda B 
= \left[\begin{array}{ccc} & & \alpha_{1} \\ 
& \alpha_{2} & \alpha_{21} \\
\alpha_{3} & \alpha_{32} & \alpha_{31} \end{array}\right]
- \lambda 
\left[\begin{array}{ccc} & & \beta_{1} \\ 
& \beta_{2} & \beta_{21} \\
\beta_{3} & \beta_{32} & \beta_{31} \end{array}\right].
\end{displaymath} 
Again we borrow the symbols $A$ and $B$ to denote the relevant submatrices.   The subpencil inherits the 
palindromic or alternating structure of the big pencil, but we will proceed at first 
as if $A$ and $B$ were unrelated.  
The poles are $\sigma_{k-1} = \alpha_{3}/\beta_{3}$, $\sigma_{k} = \alpha_{2}/\beta_{2}$, 
and $\sigma_{k+1} = \alpha_{1}/\beta_{1}$.  We want to make a unitary 
congruence transformation that swaps the positions of $\sigma_{k-1}$ and $\sigma_{k+1}$ while 
leaving $\sigma_{k}$ in the middle.    We proceed just as in the case of a $2 \times 2$ swap.   
We flip the rows and columns to get
\begin{displaymath}
FAF - \lambda FBF = 
\left[\begin{array}{ccc} 
\alpha_{31} & \alpha_{32} & \alpha_{3}  \\ \alpha_{21} & \alpha_{2} & \\ \alpha_{1} & & 
\end{array}\right]
- \lambda 
\left[\begin{array}{ccc} 
\beta_{31} & \beta_{32} & \beta_{3}  \\ \beta_{21} & \beta_{2} & \\ \beta_{1} & & 
\end{array}\right],
\end{displaymath}
which has the poles in the desired locations but the ``wrong'' triangularity.   
As in the $2 \times 2$ case, we find it convenient to work with partially flipped forms
\begin{displaymath}
AF - \lambda BF = 
\left[\begin{array}{ccc} \alpha_{1} &  &   \\
\alpha_{21}& \alpha_{2} &   \\ \alpha_{31} & \alpha_{32} & \alpha_{3} \end{array}\right]
- \lambda 
\left[\begin{array}{ccc} \beta_{1} &  &   \\
\beta_{21} & \beta_{2} &   \\ \beta_{31} & \beta_{32} & \beta_{3} \end{array}\right].
\end{displaymath}

Next we set up and solve some Sylvester equations to diagonalize the pencil.  Specifically, we will find unit
lower triangular 
\begin{displaymath}
X = \left[\begin{array}{ccc} 1 &  &  \\ x_{21} & 1 &  \\ x_{31} & x_{32} & 1 \end{array}\right]
\quad\mbox{and}\quad
Y = \left[\begin{array}{ccc} 1 & &  \\ y_{21} & 1 &  \\ y_{31} & y_{32} & 1 \end{array}\right]
\end{displaymath}
such that 
\begin{equation}\label{eq:sylvic}
(AF)X = Y(\check{A}F) \quad\mbox{and}\quad (BF)X = Y(\check{B}F),
\end{equation}
where $\check{A}$ and $\check{B}$ are anti-diagonal matrices with the same anti-diagonals as $A$ and $B$, 
respectively.  This is the $3 \times 3$ analog of (\ref{eq:sylv}).  Writing the first of these equations out in detail, we have 
\begin{displaymath}
\left[\begin{array}{ccc} \alpha_{1} &  &  \\
\alpha_{21} & \alpha_{2} &   \\ \alpha_{31} & \alpha_{32} & \alpha_{3} \end{array}\right]
\left[\begin{array}{ccc} 1 &  &  \\ x_{21} & 1 & \\ x_{31} & x_{32} & 1 \end{array}\right] =
\left[\begin{array}{ccc} 1 & & \\ y_{21} & 1 &  \\ y_{31} & y_{32} & 1 \end{array}\right]
\left[\begin{array}{ccc} \alpha_{1} &  &   \\
& \alpha_{2} &   \\ & & \alpha_{3} \end{array}\right],
\end{displaymath}
and similarly for the $B$ equation.  Altogether this is a system of six linear equations in 
the six unknowns $x_{21}$, $x_{31}$, $x_{32}$, $y_{21}$, $y_{31}$, and $y_{32}$, but fortunately 
it turns out to be three systems of two equations that are nearly independent of one another:
\begin{displaymath}
\left[\begin{array}{cc} \alpha_{1} & \alpha_{2}  \\ \beta_{1} & \beta_{2} \end{array}\right]
\left[\begin{array}{c} \phantom{-}y_{21} \\ -x_{21} \end{array}\right] = 
\left[\begin{array}{c} \alpha_{21} \\ \beta_{21} \end{array}\right],
\end{displaymath}
which has a unique solution if and only if $\alpha_{1}/\beta_{1} \neq \alpha_{2}/\beta_{2}$,
\begin{displaymath}
\left[\begin{array}{cc} \alpha_{2} & \alpha_{3} \\ \beta_{2} & \beta_{3} \end{array}\right]
\left[\begin{array}{c}  \phantom{-}y_{32} \\ - x_{32} \end{array}\right] = 
\left[\begin{array}{c} \alpha_{32} \\ \beta_{32} \end{array}\right],
\end{displaymath}
which has a unique solution if and only if $\alpha_{2}/\beta_{2} \neq \alpha_{3}/\beta_{3}$, and
\begin{displaymath}
\left[\begin{array}{cc} \alpha_{1} & \alpha_{3} \\ \beta_{1} & \beta_{3} \end{array}\right]
\left[\begin{array}{c}  \phantom{-}y_{31} \\ -x_{31} \end{array}\right] = 
\left[\begin{array}{c} \alpha_{31} + \alpha_{32}x_{21} \\ \beta_{31} + \beta_{32}x_{21} \end{array}\right],
\end{displaymath}
which has a unique solution if and only if $\alpha_{1}/\beta_{1} \neq \alpha_{3}/\beta_{3}$.
We conclude that there are unique unit lower triangular $X$ and $Y$ such that (\ref{eq:sylvic}) holds
if and only if the poles are distinct.   We can easily and stably compute $X$ and $Y$.
Rewriting (\ref{eq:sylvic}) we have the equivalence 
\begin{displaymath}
(YF)^{-1}(A - \lambda B)(FX) = F(\check{A} - \lambda\check{B})F = 
\left[\begin{array}{ccc} & & \alpha_{3}  \\ & \alpha_{2} & \\ \alpha_{1} & & \end{array}\right]
- \lambda
\left[\begin{array}{ccc} & & \beta_{3}  \\ & \beta_{2} & \\ \beta_{1} & & \end{array}\right],
\end{displaymath} 
which is the $3 \times 3$ analog of (\ref{eq:nonuneq2}).  From here on the argument is exactly the 
same as in the $2 \times 2$ case, starting from (\ref{eq:nonuneq2}), with obvious trivial modifications.  
In the end we get a congruence  $Q^{*}(A - \lambda B)Q = \hat{A} - \lambda \hat{B}$ 
that makes the desired swap.  

The computation is again straightforward.  We need a unitary $Q$ such that $FX = QR$.   The entries of $X$
can be computed by three applications of Cramer's rule.  Then the $QR$ decomposition can be computed, producing
$Q$ as a product of three core transformations.  $Q$ depends on the first two columns of $FX$, and these can
be rescaled to avoid the divisions implied by Cramer's rule without altering $Q$.  

We emphasize that this move succeeds if the three poles in question are all distinct.  In our application below,
the poles will be $\rho$, $\sigma_{n/2}$, $\tilde{\rho}$.  In the palindromic case, $\sigma_{n/2}$ is an unpaired pole 
that must lie on the unit circle, and $\tilde{\rho} = 1/\overline{\rho}$.  Thus $\rho \neq \tilde{\rho}$ if and only if $\rho$
and $\tilde{\rho}$ are not on the unit circle.  Thus, if $\rho \neq \tilde{\rho}$, then all three poles are 
automatically distinct.  The same is true in the alternating case by a similar argument involving the imaginary axis 
instead of the unit circle.  

\subsection*{Refinement of a move of type IIe}

We showed above that a move of type IIo can be refined if necessary, and the same is true of a move of type IIe.
After the move of type IIe, the resulting $3 \times 3$ pole pencil has the form 
\begin{displaymath}
A - \lambda B  = \left[\begin{array}{ccc} 
\epsilon_{31} & \epsilon_{32} & \alpha_{3} \\
\epsilon_{21} & \alpha_{2} & \alpha_{23} \\ \alpha_{1} & \alpha_{12} & \alpha_{13} \end{array}\right]
- \lambda \left[\begin{array}{ccc} 
\eta_{31} & \eta_{32} & \beta_{3} \\
\eta_{21} & \beta_{2} & \beta_{23} \\ \beta_{1} & \beta_{12} & \beta_{13} \end{array}\right].
\end{displaymath}
We have simplified the notation by leaving off the hats.  The numbers 
$\absval{\epsilon_{ij}}$ and $\absval{\eta_{ij}}$ are tiny and would have been zero except for roundoff
errors.  If they are not small enough to be ignored, we must do a refinement step.  The ``upside down'' notation
used here reflects the fact that in the analysis above, we flipped the rows of the matrices.  We could do the 
same thing here (flip the rows), but for this brief summary we will not bother.  For details see Section~\ref{sec:refine}.   

For the correction step we look for matrices
\begin{displaymath}
X = \left[\begin{array}{ccc} 1 & & \\ x_{21} & 1 & \\ x_{31} & x_{32} & 1 \end{array}\right] \quad\mbox{and}\quad
Y = \left[\begin{array}{ccc} 1 & y_{12} & y_{13} \\ & 1 & y_{23} \\ & & 1 \end{array}\right]
\end{displaymath}
that set the tiny numbers to zero:
\begin{equation}\label{eq:riccati3}
Y(A - \lambda B)X = \left[\begin{array}{ccc} 
0 & 0 & \check{\alpha}_{3} \\
0 & \check{\alpha}_{2} & \check{\alpha}_{23} \\ 
\check{\alpha}_{1} & \check{\alpha}_{12} & \check{\alpha}_{13} \end{array}\right]
- \lambda \left[\begin{array}{ccc} 
0 & 0 & \check{\beta}_{3} \\
0 & \check{\beta}_{2} & \check{\beta}_{23} 
\\ \check{\beta}_{1} & \check{\beta}_{12} & \check{\beta}_{13} \end{array}\right].
\end{equation}
Since the corrections to be made are tiny, we expect all of the the numbers $x_{ij}$ and $y_{ij}$ to be tiny.  
Writing out (\ref{eq:riccati3}) in detail we get a system of six quadratic equations in six unknowns.  The 
quadratic terms and a few others are  negligible.  Eliminating negligible terms we get six linear equations
that have a unique solution as long as the poles are distinct.  See Section~\ref{sec:refine} for details.  
In both alternating and palindromic cases, the transforming matrices $X$ and $Y$ are related by $Y = X^{*}$, 
which simplifies the situation further.  In either case, the equations are easily and stably solved.    

Once $X$ has been computed, the decomposition $X = QR$ supplies the needed unitary transforming 
matrix.  The congruence $Q^{*}(A - \lambda B)Q = \hat{A} - \lambda \hat{B}$ yields the desired refinement.  

Remark:  The moves of types IIo and IIe are special cases of constructions presented in \cite[\S\S 4.1,\,8.1]{KrScWa09}, but we
have taken a different approach here.  In \cite{KrScWa09}
the palindromic and alternating cases were considered separately.   Here we have looked at the unstructured case first
and shown how to do the swap using a unitary equivalence.   Then we have shown that in both of our structured cases
the equivalence is actually a congruence.     
We have also added a refinement step, which was not contemplated in \cite{KrScWa09}, to make the algorithm 
more robust.  

\section{Stability}

We can construct a variety of algorithms from the moves.   If each move is backward stable, then any algorithm
built from moves must also be backward stable.  We therefore take a moment to consider this question.  
Standard backward error analysis \cite{Wil65} shows that moves of type I are backward stable.  
If the moves of type II are implemented as shown in \cite{CaMaVaWa19}, they never fail and are 
always backward stable.  

Moves of types IIo and IIe require the solution of small linear systems that are nonsingular if and 
only if the poles involved in the swap are distinct.  As we will see below, there are other good reasons
(involving convergence rates) for keeping these poles distinct and preferably far apart.  Assuming this 
is done, we can expect these moves to be stable.  There is a natural stability test associated with these moves:
check that the numbers that are supposed to be zero really are (almost) zero.  For the 
event that they are not, we have described a refinement step that can be used to make them smaller.  
The refinement can be repeated if necessary, though this should be vary rare.  Because of the refinement
step, we can say for sure that moves of types IIo and IIe are backward stable, provided that we do not attempt
to swap two equal poles.  

\section{Building an algorithm using the moves}\label{sec:algorithm}

First suppose our pair $(A,B)$ has odd dimension $n$, and its poles are 
$\sigma_{1}$, $\sigma_{2}$, \ldots, $\sigma_{m}$, $\tilde{\sigma}_{m}$, \ldots, $\tilde{\sigma}_{2}$, $\tilde{\sigma}_{1}$, 
where $m = (n-1)/2$, and 
$\tilde{\sigma}_{i} = 1/\overline{\sigma}_{i}$ (resp.\ $-\overline{\sigma}_{i}$) in the palindromic (resp.\ alternating) case.  
One iteration of the most basic algorithm would proceed as 
follows.  First a shift $\rho$ is chosen and inserted in place of $\sigma_{1}$ by a move of type I.  This move also inserts 
a shift $\tilde{\rho}$ in place of $\tilde{\sigma}_{1}$ at the other end.   A simple choice of $\rho$ would be the 
\emph{Rayleigh quotient shift} $\rho = a_{1,n}/b_{1,n}$.  Then a move of type II is used to interchange
$\rho$ with $\sigma_{2}$ and $\tilde{\rho}$ with $\tilde{\sigma}_{2}$.  Then another move of type II is used to interchange
$\rho$ with $\sigma_{3}$, and so on.  After $m-1$ such moves, the poles will be 
$\sigma_{2}$, \ldots, $\sigma_{m}$, $\rho$, $\tilde{\rho}$, $\tilde{\sigma}_{m}$, \ldots, $\tilde{\sigma}_{2}$, with $\rho$ 
and $\tilde{\rho}$ side by side in the middle.  
Then a move of type IIo can be used to swap them, provided that $\rho \neq \tilde{\rho}$.   Then additional moves of type II are used to push $\rho$ and $\tilde{\rho}$ further along, that is, $\rho$ is swapped with $\tilde{\sigma}_{m}$, 
then $\tilde{\sigma}_{m-1}$, and so on, while $\tilde{\rho}$ is swapped with $\sigma_{m}$, $\sigma_{m-1}$, etc.  
Once the shifts arrive at the edge of the pencil, they can be removed by a move of type I, which would replace them by new poles, which could be the original poles $\sigma_{1}$, $\tilde{\sigma_{1}}$, or they could be different.  This completes
one iteration of the basic algorithm.

The case of even dimension is the same, except that there is an extra unpaired pole $\sigma_{n/2}$ in the middle.  
Type II operations push the shifts $\rho$ and $\tilde{\rho}$ toward the middle, as in the odd case, until the configuration
of poles is $\sigma_{2}$, \ldots, $\rho$, $\sigma_{n/2}$, $\tilde{\rho}$, \ldots, $\tilde{\sigma_{2}}$.   Then a move of 
type IIe is used to swap $\rho$ and $\tilde{\rho}$  while leaving $\sigma_{n/2}$ fixed.   This can be done 
provided that $\rho \neq \tilde{\rho}$.  Once this exchange has been made, the iteration is completed in the same way
as in the odd case.  

Repeated iterations of the basic algorithm with good choices of shifts $\rho$ and $\tilde{\rho}$ will cause the pencil
to tend toward triangular form, exposing the eigenvalues on the anti-diagonal.  Not all eigenvalues appear at once.  
Good shifts can (usually) cause $a_{n-1,1} \to 0$ and $b_{n-1,1} \to 0$ in just a few iterations; the convergence rate is
typically quadratic.  By symmetry we must also have $a_{1,n-1} \to 0$ and $b_{1,n-1} \to 0$ at the same rate.  
This exposes a pair of eigenvalues  $\lambda = a_{n,1}/b_{n,1}$ and $\tilde{\lambda} = a_{1,n}/b_{1,n}$ at the ends.  
Then the problem can be deflated to size $n-2$, and we can go after the next pair of eigenvalues, and so on.  All of
this is a consequence of the convergence theorem stated below.   We will not present a detailed explanation 
because the arguments are the same as in the unstructured case.   

\subsection*{A convergence theorem}

The mechanism that drives all variants of Francis's algorithm, including the QZ algorithm, is nested subspace 
iteration with changes of coordinate system.  See 
\cite[p.~431]{Wat10}, \cite[p.~399]{Wat11}, or \cite[Theorem~2.2.3]{AuMaRoVaWa18}.  
This is also true of our basic algorithm sketched above.  We just need to take a few lines to
set the scene.  We make the (generically valid) assumption that none of the poles or shifts is exactly an eigenvalue
of the pencil $A - \lambda B$.   We continue to cover the alternating and palindromic cases simultaneously.  
Each iteration begins with the choice of a shift $\rho$ and a companion shift $\tilde{\rho}$ 
($= 1/\overline{\rho}$ in the palindromic case and $-\overline{\rho}$ in the alternating case).   
The result of the iteration is a new structured
anti-Hessenberg pencil $\hat{A} - \lambda \hat{B}$ satisfying 
\begin{equation}\label{eq:palcoord}
\hat{A} - \lambda \hat{B} = Q^{*}(A - \lambda B)Q.
\end{equation}
We need to define two nested sequences of subspaces.  
For $k=1$, \ldots, $n$, define 
\begin{displaymath}
\eee_{k} = \spn{e_{1}, \ldots, e_{k}},
\end{displaymath}
where $e_{1}$, \ldots, $e_{n}$ are the standard basis vectors.  Then define 
\begin{displaymath}
\cue_{k} = Q\eee_{k},
\end{displaymath}
the space spanned by the first $k$ columns of $Q$.

\begin{theorem}\label{thm:palspace}
A single step of the basic algorithm described above with shift $\rho$ effects nested subspace iterations  
\begin{displaymath} 
\cue_{k} = (A - \tilde{\rho}B)^{-1}(A - \rho B)\eee_{k}, \qquad  \qquad k = 1,\ldots,n-1.
\end{displaymath}
The change of coordinate system (\ref{eq:palcoord}) transforms $\cue_{k}$ back to $\eee_{k}$.
\end{theorem}

This theorem makes no mention of convergence, but we call it a \emph{convergence theorem} anyway.  
This result and ones like it can be used together with the convergence theory of 
subspace iteration to draw conclusions about
the convergence of the algorithm, as explained in \cite{Wat07,Wat10,Wat11} and elsewhere.  

\begin{proof}
We sketch the proof, relying on Theorem~5.2 of \cite{CaMaVaWa19}.  That theorem applies to upper-Hessenberg 
pencils, so we can apply it to the flipped pencil $FA - \lambda FB$.  Notice that 
$(FA - \tilde{\rho}FB)^{-1}(FA - \rho FB) = (A - \tilde{\rho}B)^{-1}(A - \rho B)$.  

First suppose $n$ is odd and $k < (n-1)/2$.  According to  (the ``$Z$'' part of) 
Theorem~5.2 of \cite{CaMaVaWa19}, the action on $\eee_{k}$
depends on the two moves of type II that take place at the ``$k$th position'', by which we mean the spot originally 
occupied by poles $\sigma_{k}$ and $\sigma_{k+1}$.   The basic algorithm inserts the shift $\rho$, swapping it with
$\sigma_{1}$, \ldots, $\sigma_{k}$.   The first swap at the $k$th position is an exchange of $\rho$ with $\sigma_{k+1}$,
which (see \cite[Theorem~5.2]{CaMaVaWa19}) generates a factor $(z-\rho)/(z - \sigma_{k+1})$.  The only other swap 
at the $k$th position happens when $\sigma_{k+1}$ is swapped with $\tilde{\rho}$, moving $\sigma_{k+1}$ back to
its original position.  This introduces a factor $(z - \sigma_{k+1})/(z - \tilde{\rho})$.  The action on $\eee_{k}$ is determined 
by the product of the factors, which is $r(z) = (z - \rho)/(z - \tilde{\rho})$.  Specifically, $\eee_{k}$ is transformed to 
$\cue_{k} = r(B^{-1}A)\eee_{k}$.  Since $r(B^{-1}A) = (B^{-1}A - \tilde{\rho}I)^{-1}(B^{-1}A - \rho I) = 
(A - \tilde{\rho}B)^{-1}(A - \rho B)$, we have $\cue_{k} = (A - \tilde{\rho}B)^{-1}(A - \rho B)\eee_{k}$, as claimed.

The case $k > (n-1)/2$ is the same, except that the order of the swaps is reversed.  This does not change the outcome.  

The case $k = (n-1)/2$ is different because this is the spot in the middle where only a single swap 
takes place, interchanging the shifts $\rho$ and $\tilde{\rho}$.  
Using \cite[Theorem~5.2]{CaMaVaWa19} again, we see that we 
just have a single factor $(z - \rho)/(z - \tilde{\rho})$, so the result is the same as in the other cases.    
This case can also be deduced directly from \cite[Theorem~4.3]{CaMaVaWa19}, which is a precursor of 
\cite[Theorem~5.2]{CaMaVaWa19}.

In the case of even $n$, the argument is the same as in the odd case if $k < n/2 -1$ or $k > n/2$.  The only question
mark is in the cases $k = n/2-1$ and $k= n/2$, which are affect by (and only by) the move of type IIe in the middle.
Since this move, which affects three adjacent poles instead of two, is different from a standard move of type II, we 
must check what happens here.  We leave it to the reader to verify that the proof of 
\cite[Theorem~4.3]{CaMaVaWa19}, 
which applies to moves of type II, is also valid for moves of type IIe.  Once this has been checked, the proof is 
complete.   \hfill\end{proof}

In our discussion of moves of type IIo and IIe we had to make the assumption $\rho \neq \tilde{\rho}$ in order to 
ensure that the moves are possible.  Theorem~\ref{thm:palspace} gives us another reason for this assumption:  
If $\rho = \tilde{\rho}$, we have $(A - \tilde{\rho}B)^{-1}(A - \rho B)=I$, and the iteration goes nowhere.  The action of the shift $\rho$
traveling in one direction is exactly cancelled by the shift $\tilde{\rho}$ traveling in the opposite direction.  

In the alternating case the requirement $\rho \neq \tilde{\rho}$ means that the shifts should not lie on the imaginary axis.
It follows that this method will only be useful for finding eigenvalues that are not purely imaginary.  This is not necessarily 
a weakness.  For the continuous-time optimal control problems mentioned in the introduction \cite{Meh91}, mild assumptions 
guarantee that no eigenvalues are on the imaginary axis;  exactly half are in the open left half plane and half are in the
open right half plane.  Our basic algorithm will have no problem computing all of these eigenvalues;  they will be 
extracted in $(\lambda,-\overline{\lambda})$ pairs.

In the palindromic case the requirement $\rho \neq \tilde{\rho}$ means that the shifts should not lie on the unit 
circle, so this method will only be useful for finding eigenvalues that have modulus different from 1.  For the 
discrete-time control problems mentioned above, mild assumptions guarantee that no eigenvalues lie on the 
unit circle; exactly half are inside and half are outside.  Our algorithm will easily compute all of these eigenvalues, 
and they will be extracted in $(\lambda,1/\overline{\lambda})$ pairs.

The bulge-chasing algorithm in \cite{KrScWa09} has the same ($\rho \neq \tilde{\rho}$) restriction.  In fact we can
show that our basic algorithm is a generalization of the algorithm in 
\cite{KrScWa09}.\footnote{In \cite{KrScWa09}  we considered multi-shift bulge-chasing algorithms of arbitrary degree.
Here we are considering only a single-shift algorithm, and this generalizes  the single-shift version of the algorithm in
\cite{KrScWa09}.}

Let's take a look at this algorithm, beginning with the palindromic case.  
The pencil is  $A - \lambda A^{*}$ with $A$ assumed to be \emph{anti-Hessenberg}, 
i.e.\
\begin{displaymath}
A = 
\parbox{3.2cm}{
\begin{tikzpicture}[scale=1.66,y=-1cm]
\draw (-.15,-.1) -- (-.2,-.1) -- (-.2,1.5) -- (-.15,1.5);
\draw (1.55,-.1) -- (1.6,-.1) -- (1.6,1.5) -- (1.55,1.5);
\foreach \j in {0,...,7}{
   \foreach \i in {\j,...,7}{\node at (\i/5,1.4-\j/5)
     [align=center,scale=1.0]{$\times$};}}
\foreach \j in {0,...,6}{\node at (\j/5,1.4-\j/5-.2)
     [align=center,scale=1.0]{$\times$};}
\end{tikzpicture}
}.
\end{displaymath}
The algorithm presented in \cite{KrScWa09} 
begins with a reduction step that introduces some zeros above the anti-diagonal, 
transforming $A$ to a partially anti-triangular form 
\begin{equation}\label{eq:partanttri}
\parbox{3.2cm}{
\begin{tikzpicture}[scale=1.66,y=-1cm]
\draw (-.15,-.1) -- (-.2,-.1) -- (-.2,1.5) -- (-.15,1.5);
\draw (1.55,-.1) -- (1.6,-.1) -- (1.6,1.5) -- (1.55,1.5);
\foreach \j in {0,...,7}{
   \foreach \i in {\j,...,7}{\node at (\i/5,1.4-\j/5)
     [align=center,scale=1.0]{$\times$};}}
\foreach \j in {0,...,3}{\node at (\j/5,1.4-\j/5-.2)
     [align=center,scale=1.0]{$\times$};}
\end{tikzpicture}
}.
\end{equation}
A bulge-chasing algorithm is then applied to this partially reduced form to expose the  
eigenvalues $\{\lambda,\, \overline{\lambda}^{-1}\}$ satisfying $\absval{\lambda} \neq 1$
in pairs.   

From our current vantage point we can see that the preliminary reduction step is just 
a process of introducing zero and infinite poles into the pencil $A - \lambda A^{*}$ by moves
of types I and II.   This step is necessary for the bulge-chasing algorithm in \cite{KrScWa09},
but it is not needed for the pole swapping algorithm discussed in this paper;
we can go to work right away on the anti-Hessenberg pencil $A - \lambda A^{*}$.  

The algorithm in \cite{KrScWa09} for the alternating case is a little bit different.  Starting with
an alternating pencil $A - \lambda B$ in anti-Hessenberg form, this algorithm requires a preliminary 
step to transform $B$ to anti-triangular form, leaving $A$ anti-Hessenberg.  Then a bulge-chasing
algorithm is applied.   The modified pencil has all poles equal to infinity.  
Again, from our new viewpoint, we can see that the preliminary reduction
is nothing but a process of introducing infinite poles by moves of types I and II.  This is necessary
for the bulge-chasing algorithm in \cite{KrScWa09}, but it is not needed for our pole-swapping algorithm.

In both cases the flop count for this additional reduction is $O(n^{3})$  ($O(n^{2})$ moves at $O(n)$ flops per move),
so the cost is not insignificant.

For a single-shift iteration of the algorithm in \cite{KrScWa09} we have a theorem just like Theorem~\ref{thm:palspace}.
If we introduce a shift $\rho$ at one end, we automatically introduce a complementary shift $\tilde{\rho}$ 
( $= 1/\overline{\rho}$ or $-\overline{\rho}$) as always, at the other end.  
The setup is the same as for Theorem~\ref{thm:palspace}, and the new result looks like this:

\begin{theorem}\label{thm:palspace2}
A single step of the algorithm in \cite{KrScWa09} with shift $\rho$ effects nested subspace iterations  
\begin{displaymath} 
\cue_{k} = (A - \tilde{\rho}B)^{-1}(A - \rho B)\eee_{k}, \qquad  \qquad k = 1,\ldots,n-1.
\end{displaymath}
The change of coordinate system (\ref{eq:palcoord}) transforms $\cue_{k}$ back to $\eee_{k}$.
\end{theorem}

This is an immediate consequence of \cite[Theorem~7.3.1]{Wat07}.   Note that Theorem~\ref{thm:palspace2} is identical to
Theorem~\ref{thm:palspace}.   The only difference is that the bulge-chasing 
algorithm requires a special Hessenberg-triangular form, as described immediately above.   Since the action of the 
algorithms is the same, we deduce that the pole-swapping algorithm is a generalization of the bulge-chasing 
algorithm.    

Remark:  For clarity we have focused on the most basic possible pole-swapping algorithm for the palindromic and alternating 
problems.  One can consider variants that introduce multiple shifts and draw the same conclusions.  
For more ideas see \cite{CaMaVaWa19}, which also contains a more explicit demonstration of the 
connection between bulge-chasing and pole-swapping algorithms.

\section{Justification of the refinement step}\label{sec:refine}

In Section~\ref{sec:operations} we briefly described refinement procedures that can be applied (occasionally)
after moves of types IIo and IIe.  Here we provide complete justifications for those procedures.   
The algorithm in \cite{KrScWa09} also has a move in the middle to which a refinement step could be applied.
In that paper it was acknowledged that a failure might
occasionally occur, but it was reported that in the course of the various tests, no failures were observed.   
However, since a failure might occur at any time, it would make sense to add the refinement step to the algorithm
of \cite{KrScWa09} to make it more robust.  In order to accommodate its use in the context of \cite{KrScWa09},
we have made our discussion of the refinement procedure more general than is strictly required for this paper.

First we consider the case of no middle pole(s), as in a move of type IIo.  Suppose we have just swapped $m$ poles 
with $m$ other poles.  After the swap we have a bulge pencil
\begin{displaymath}
\left[\begin{array}{cc} 
E_{11} & A_{12} \\ A_{21} & A_{22}
\end{array}\right] - \lambda
\left[\begin{array}{cc} 
G_{11} & B_{12} \\ B_{21} & B_{22}
\end{array}\right],
\end{displaymath}
where the submatrices are $m \times m$.  
The matrices $E_{11}$ and $G_{11}$ would be zero if there were no roundoff errors,  
so $\norm{E_{11}} \ll \norm{A}$ and $\norm{G_{11}} \ll \norm{B}$.   We assume that the eigenvalues (the poles)
of the subpencil $A_{21} - \lambda B_{21}$ are disjoint from those of  $A_{12} - \lambda B_{12}$.   In the case
of a move of type IIo we have $m=1$, but we are now allowing larger $m$ to take into account the scenario of
\cite{KrScWa09}.  

If $\norm{E_{11}}$ and $\norm{G_{11}}$ are not small enough, we must do a refinement step.
To this end we seek $X$ and $Y$ such that 
\begin{equation}\label{eq:genric1}
\begin{array}{c}
\left[\begin{array}{cc} I & Y \\ & I \end{array}\right]
\left[\begin{array}{cc} E_{11} & A_{12} \\  A_{21} & A_{22} \end{array}\right]
\left[\begin{array}{cc} I & \\ X & I \end{array}\right] = 
\left[\begin{array}{cc}  0 & \check{A}_{12}  \\
\check{A}_{21} & \check{A}_{22} \end{array}\right],  \\ \\
\left[\begin{array}{cc} I & Y \\ & I \end{array}\right]
\left[\begin{array}{cc} G_{11} & B_{12} \\ B_{21} & B_{22} \end{array}\right]
\left[\begin{array}{cc} I & \\ X & I \end{array}\right] = 
\left[\begin{array}{cc}  0 & \check{B}_{12} \\ \check{B}_{21} & \check{B}_{22}
\end{array}\right]. 
\end{array}  
\end{equation}
By straightforward computation we find that (\ref{eq:genric1}) holds if and only if the algebraic Riccati equations 
\begin{equation} \label{eq:genric2}
\begin{array}{ccc}
A_{12}X + YA_{21} + E_{11} + YA_{22}X & = & 0   \\
B_{12}X + YB_{21} + G_{11} + YB_{22}X & = & 0 
\end{array} 
\end{equation}
hold.  
Since $\norm{E_{11}}$ and $\norm{G_{11}}$ are tiny, we expect that the corrections 
$X$ and $Y$ will be tiny as well.  Therefore the quadratic terms $YA_{22}X$ and $YB_{22}X$ in (\ref{eq:genric2}) 
should be negligible, and (\ref{eq:genric2}) should be well approximated by the Sylvester equations 
\begin{equation} \label{eq:gensyl1}
\begin{array}{ccc}
A_{12}X + YA_{21} + E_{11}  & = & 0   \\
B_{12}X + YB_{21} + G_{11} & = & 0. 
\end{array} 
\end{equation}
These linear equations have a unique solution $(X,Y)$ if and only if the eigenvalues of the pencil 
$A_{12} - \lambda B_{12}$ are disjoint
from those of $A_{21} - \lambda B_{21}$ \cite[Theorem~6.6.8]{Wat07}.  Under this assumption one can prove 
by standard arguments  \cite{Ste72,Ste73}, \cite[\S~2.7]{Wat07} (using the contraction mapping principle)
that the Riccati equations (\ref{eq:genric2}) have 
a unique small solution $(X,Y)$ if $\norm{E_{11}}$ and $\norm{G_{11}}$ are sufficiently small.  By this we mean
that, although (\ref{eq:genric2}) typically has many solutions,  
there is exactly one for which $\norm{X}$ and $\norm{Y}$ are small, and this is the solution that is of interest to us. 

So far we have ignored the special structure of the pencil.  Now let's see what we can say in the alternating 
case, for which $A_{12}^{*} = A_{21}$, $A_{22}^{*} = A_{22}$, $E_{11}^{*} = E_{11}$, $B_{12}^{*} = -B_{21}$, 
$B_{22}^{*} = -B_{22}$, and $G_{11}^{*} = -G_{11}$.  If we make these substitutions in (\ref{eq:genric2}) and then 
take conjugate transposes, we get 
\begin{equation}\label{eq:genric3}
\begin{array}{ccc} 
A_{12}Y^{*} + X^{*}A_{21} + E_{11} +X^{*}A_{22}Y^{*} = 0 \\
B_{12}Y^{*} + X^{*}B_{21} + G_{11} +X^{*}B_{22}Y^{*} = 0. 
\end{array}
\end{equation}

Now consider the palindromic case, for which $A_{12}^{*}=B_{21}$, $A_{21}^{*} = B_{12}$, $A_{22}^{*} = B_{22}$, 
and $E_{11}^{*} = G_{11}$.  Making these substitutions in (\ref{eq:genric2}) and then taking conjugate transposes,
we find again that we get (\ref{eq:genric3}).  

From (\ref{eq:genric3}) we deduce that, in both of our structured cases, 
$(Y^{*},X^{*})$ is a solution of (\ref{eq:genric2}) if and only if $(X,Y)$ is. 
Since (\ref{eq:genric2}) has a unique small-norm solution, we deduce that $(Y^{*},X^{*}) = (X,Y)$, or briefly
$Y = X^{*}$.   If we now make this substitution  into (\ref{eq:genric1}), we see that the transformation is a congruence,
which is exactly what is needed for the preservation of the two structures.  

If we solve (\ref{eq:genric2}) for $X$, we can carry out the congruence (\ref{eq:genric1}) to make the desired correction.
However, (\ref{eq:genric1}) has the shortcoming that it is not unitary.  To get a unitary congruence  
that has the same effect, we can do a QR decomposition
\begin{equation}\label{eq:qr}
\left[\begin{array}{cc} I & \\ X & I\end{array}\right] = QR = 
\left[\begin{array}{cc} Q_{11} & Q_{12} \\ Q_{21} & Q_{22} \end{array}\right]
\left[\begin{array}{cc} R_{11} & R_{12} \\ & R_{22} \end{array}\right],
\end{equation}
where $Q$ is unitary and $R$ is upper triangular, and do a congruence with $Q$.  Substituting $QR$ into 
(\ref{eq:genric1}) in two places and inverting the triangular matrices $R^{*}$ and $R$, we obtain 
\begin{displaymath}
Q^{*}
\left[\begin{array}{cc} E_{11} & A_{12} \\  A_{21} & A_{22} \end{array}\right]
Q = 
R^{-*}\left[\begin{array}{cc}  0 & \check{A}_{12}  \\
\check{A}_{21} & \check{A}_{22} \end{array}\right]R^{-1}  = 
\left[\begin{array}{cc} 0 & \hat{A}_{12} \\ \hat{A}_{21} & \hat{A}_{22} \end{array}\right],
\end{displaymath}
where 
\begin{displaymath}
\hat{A}_{12} = R_{11}^{-*}\check{A}_{12}R_{22}^{-1} \quad\mbox{and}\quad 
\hat{A}_{21} = R_{22}^{-*}\check{A}_{21}R_{11}^{-1}.
\end{displaymath}
We have displayed only the ``$A$'' equations, but the ``$B$'' equations are the same.

We now have all of the ingredients we need for an update.  
The discussion so far suggests that we will compute $X$ by solving
the Riccati equations (\ref{eq:genric2}), but in fact we will not.  Instead we will solve the Sylvester equations
(\ref{eq:gensyl1}) to get an excellent approximation.  This amounts to one step of Newton's method applied to 
(\ref{eq:genric2}) using initial guess  $X^{(0)} = Y^{(0)} = 0$.  Notice that the symmetry argument that we applied 
to (\ref{eq:genric2}) above also applies to the Sylvester equation:  In both the alternating and palindromic cases, 
if $(X,Y)$ is the unique solution of (\ref{eq:gensyl1}), then so is $(Y^{*},X^{*})$, so $Y = X^{*}$.  

When we solve (\ref{eq:gensyl1}) we can take symmetry into account.  In the palindromic case, the two 
matrix equations of (\ref{eq:gensyl1}) are equivalent, so we just have to solve 
\begin{displaymath}
A_{12}X + X^{*}A_{21} + E_{11} = 0
\end{displaymath}
for $X$.  Separating real and imaginary parts, we can write this as a system of $2m^{2}$ real linear equations
in $2m^{2}$ unknowns, the real and imaginary parts of $X$.  This can be solved stably by conventional means.
In the alternating case the $A$ equation is symmetric and the $B$ equation is skew symmetric.   Taking these 
symmetries into account, we again get a system of $2m^{2}$ real equations in $2m^{2}$ real unknowns, which
can be solved stably by conventional means.

We can now summarize our refinement step:  Solve the Sylvester equations for $X$, taking the symmetry into
account.  Then perform the QR decomposition (\ref{eq:qr}), and use the resulting $Q$ to effect a 
unitary congruence transform
\begin{displaymath}
Q^{*}
\left[\begin{array}{cc} E_{11}  & A_{12} \\
A_{21}  & A_{22}  \end{array}\right] Q = 
\left[\begin{array}{cc} \hat{E}_{11}  & \hat{A}_{12}  \\
\hat{A}_{21}  & \hat{A}_{22}  \end{array}\right], \quad
Q^{*}
\left[\begin{array}{cc} G_{11} & B_{12} \\
B_{21}  & B_{22} \end{array}\right] Q = 
\left[\begin{array}{cc}  \hat{G}_{11} & \hat{B}_{12} \\
\hat{B}_{21} & \hat{B}_{22} \end{array}\right].
\end{displaymath}
We can exploit symmetry in this step as well.  
The matrices $\hat{E}_{11}$ and $\hat{G}_{11}$ are nonzero because $X$ does not exactly satisfy 
the Riccati equations (\ref{eq:genric2}), but the quadratic convergence of Newton's method guarantees
that $\norm{\hat{E}_{11}} \ll \norm{E_{11}}$ and $\norm{\hat{G}_{11}} \ll \norm{G_{11}}$.  
Thus $\hat{E}_{11}$ and $\hat{G}_{11}$ will be small enough that they can be set to zero without compromising stability.   In the extremely rare 
event that they are not small enough, the refinement step can be repeated.  

Now we consider the refinement step in the case when we have three blocks, as in a move of type IIe.  
After the move we have a pole pencil
\begin{displaymath}
\left[\begin{array}{ccc} E_{11} & E_{12} & A_{13} \\ E_{21} & A_{22} & A_{23} \\ 
A_{31} & A_{32} & A_{33} \end{array}\right]  - \lambda 
\left[\begin{array}{ccc} G_{11} & G_{12} & B_{13} \\ G_{21} & B_{22} & B_{23} \\ 
B_{31} & B_{32} & B_{33} \end{array}\right].
\end{displaymath}
In a move of type IIe, all of the submatrices are $1 \times 1$.  Here we allow them to be of any size.
Let's say the matrices $A_{31}$, $A_{13}$, etc.\  are $m \times m$, and the matrices $A_{22}$ etc.\ 
are $k \times k$.  If we take $k=0$, this reduces to the case that we have just covered.  

The submatrices $E_{ij}$ and $G_{ij}$ would be zero except for roundoff errors.  If they are small enough,
we can set them to zero and proceed.  Otherwise we must do a refinement step, and to this end we seek
$X_{ij}$ and $Y_{ij}$ such that 
\begin{displaymath}
\left[\begin{array}{ccc} I & Y_{12} & Y_{13} \\ & I & Y_{23} \\ & & I\end{array}\right]
\left[\begin{array}{ccc} E_{11} & E_{12} & A_{13} \\ E_{21} & A_{22} & A_{23} \\ 
A_{31} & A_{32} & A_{33} \end{array}\right] 
\left[\begin{array}{ccc} I & &  \\ X_{21} & I & \\ X_{31} & X_{32} & I \end{array}\right] =
\left[\begin{array}{ccc} &  & \check{A}_{13} \\ & \check{A}_{22} & \check{A}_{23} \\ 
\check{A}_{31} & \check{A}_{32} & \check{A}_{33} \end{array}\right]
\end{displaymath}
and 
\begin{displaymath}
\left[\begin{array}{ccc} I & Y_{12} & Y_{13} \\ & I & Y_{23} \\ & & I\end{array}\right]
\left[\begin{array}{ccc} G_{11} & G_{12} & B_{13} \\ G_{21} & B_{22} & B_{23} \\ 
B_{31} & B_{32} & B_{33} \end{array}\right] 
\left[\begin{array}{ccc} I & &  \\ X_{21} & I & \\ X_{31} & X_{32} & I \end{array}\right] =
\left[\begin{array}{ccc} &  & \check{B}_{13} \\ & \check{B}_{22} & \check{B}_{23} \\ 
\check{B}_{31} & \check{B}_{32} & \check{B}_{33} \end{array}\right].
\end{displaymath}
This results in a system of six Riccati equations, three from $A$ and three from $B$:
\begin{displaymath}
\begin{array}{c}
\begin{array}{c}
A_{13}X_{31} + Y_{13}A_{31} + E_{11} + Y_{13}A_{33}X_{31} + Y_{12}E_{21} + E_{12}X_{21}  \hspace{40pt} \\  
\hspace{93pt} +Y_{12}A_{22}X_{21} + Y_{13}A_{32}X_{21} +Y_{12}A_{23}X_{31} = 0,
\end{array} \\ \\
\begin{array}{c}
B_{13}X_{31} + Y_{13}B_{31} + G_{11} + Y_{13}B_{33}X_{31} + Y_{12}G_{21} + G_{12}X_{21}  \hspace{40pt} \\  
\hspace{93pt} +Y_{12}B_{22}X_{21} + Y_{13}B_{32}X_{21} +Y_{12}B_{23}X_{31} = 0,
\end{array}\end{array}
\end{displaymath}
\begin{displaymath}
\begin{array}{c}
A_{13}X_{32} + Y_{12}A_{22} + E_{12} + Y_{12}A_{23}X_{32} + Y_{13}A_{32} + Y_{13}A_{33}X_{32} = 0, \\ \\
B_{13}X_{32} + Y_{12}B_{22} + G_{12} + Y_{12}B_{23}X_{32} + Y_{13}B_{32} + Y_{13}B_{33}X_{32} = 0, 
\end{array}
\end{displaymath}
\begin{displaymath}
\begin{array}{c}
A_{22}X_{21} + Y_{23}A_{31} + E_{21} + Y_{23}A_{32}X_{21} + A_{23}X_{31} +  Y_{23}A_{33}X_{31} = 0, \\ \\ 
B_{22}X_{21} + Y_{23}B_{31} + G_{21} + Y_{23}B_{32}X_{21} + B_{23}X_{31} +  Y_{23}B_{33}X_{31} = 0.
\end{array}
\end{displaymath}
We do not propose to compute the exact solution of these equations.  Instead we will obtain an excellent approximation
by solving the Sylvester equations that one gets by discarding all negligible terms.  For example, in the first equation
above, only the first three terms are non-negligible.  We obtain 
\begin{equation}\label{eq:gensyl11}
\begin{array}{c}
A_{13}X_{31} + Y_{13}A_{31} + E_{11} = 0, \\ \\
B_{13}X_{31} + Y_{13}B_{31} + G_{11} = 0, 
\end{array}
\end{equation}
\begin{equation}\label{eq:gensyl12}
\begin{array}{c}
A_{13}X_{32} + Y_{12}A_{22} + E_{12} + Y_{13}A_{32}  = 0, \\ \\
B_{13}X_{32} + Y_{12}B_{22} + G_{12} + Y_{13}B_{32}  = 0, 
\end{array}
\end{equation}
\begin{equation}\label{eq:gensyl21}
\begin{array}{c}
A_{22}X_{21} + Y_{23}A_{31} + E_{21} + A_{23}X_{31}  = 0, \\ \\ 
B_{22}X_{21} + Y_{23}B_{31} + G_{21} + B_{23}X_{31}  = 0.
\end{array}
\end{equation}
Equations (\ref{eq:gensyl11}) are independent of the others.  They have a unique solution $(X_{31},Y_{13})$ if and
only if the pencils $A_{13} - \lambda B_{13}$ and $A_{31} - \lambda B_{31}$ have disjoint spectra.  Once we have 
solved these equations, we can substitute $Y_{13}$ into (\ref{eq:gensyl12}) and $X_{31}$ into (\ref{eq:gensyl21}).  
Equations (\ref{eq:gensyl12}) have a unique solution $(X_{32},Y_{12})$ if and only if the pencils 
$A_{13} - \lambda B_{13}$ and $A_{22} - \lambda B_{22}$ have disjoint spectra.  Similarly (\ref{eq:gensyl21}) have
a unique solution $(X_{21},Y_{23})$ if and only if 
the spectra of $A_{22} - \lambda B_{22}$ and $A_{31} - \lambda B_{31}$ are disjoint.

In both the alternating and palindromic cases, one can show that $Y_{12} = X_{21}^{*}$, $Y_{13} = X_{31}^{*}$, 
and $Y_{23} = X_{32}^{*}$.   The routine but tedious proof is left for the reader.   When we solve the Sylvester equations
in practice, we take these symmetries into account.  For example, in the palindromic case, the $A$ and $B$ equations
are equivalent, so we only have to solve the $A$ equations;  (\ref{eq:gensyl11}) reduces to 
\begin{displaymath}
A_{13}X_{31} + X_{31}^{*}A_{31} + E_{11} = 0, 
\end{displaymath}
and (\ref{eq:gensyl12}) and (\ref{eq:gensyl21}) together reduce to 
\begin{displaymath}
\begin{array}{c}
A_{13}X_{32} + X_{21}^{*}A_{22} + E_{12} + X_{31}^{*}A_{32}  = 0, \\ \\
A_{22}X_{21} + X_{32}^{*}A_{31} + E_{21} + A_{23}X_{31}  = 0,  
\end{array}
\end{displaymath}
which can be solved simultaneously.  

Once we have computed $X$, we perform a decomposition
\begin{displaymath}
\left[\begin{array}{ccc} I & & \\ X_{21} & I & \\ X_{31} & X_{32} & I \end{array}\right]  = QR,
\end{displaymath}
and use $Q$ to do a unitary congruence 
\begin{displaymath}
Q^{*} \left[\begin{array}{ccc} E_{11} & E_{12} & A_{13}  \\
E_{21} & A_{22} & A_{23} \\ A_{31} & A_{32} & A_{33}\end{array}\right] Q =
\left[\begin{array}{ccc} \hat{E}_{11} & \hat{E}_{12} & \hat{A}_{13}  \\
\hat{E}_{21} & \hat{A}_{22} & \hat{A}_{23} \\ \hat{A}_{31} & \hat{A}_{32} & \hat{A}_{33}
\end{array}\right],
\end{displaymath}
and similarly for $B$.  Because the correction amounts to one step of Newton's method, and the errors were 
tiny to begin with, we will have $\norm{\hat{E}_{11}} \ll \norm{E_{11}}$ and so on, and we can safely set the 
new errors to zero.  In the highly unlikely event that the errors are still too big, we can repeat the refinement step.

\section{Numerical experiments}

We wrote and tested MATLAB code for the palindromic case.  We tested the most basic pole-swapping algorithm with a 
single shift, as described in Section~\ref{sec:algorithm}.  We used Wilkinson shifts \cite{Wat10}. 
For moves of types IIo and IIe we set a 
tolerance of $10\epsilon\norm{M}_{F}$, where $\epsilon$ is the machine epsilon, which is approximately 
$2.22 \times 10^{-16}$, and $M$ is the small submatrix in which the swap occurs.  If all of the numbers that are
supposed to be zero are below this tolerance, no refinement is needed.  Our MATLAB code is publicly available at 
\textit{github.com/thijssteel/palindromic-rqz}.

For our first test we considered palindromic pencils $A - \lambda A^{*}$, where $A$ is a random anti-Hessenberg 
matrix.  Each of the nonzero entries is generated as $2a + bi$, where $a$ and $b$ are drawn from a standard normal 
distribution.   We computed the eigenvalues using our pole-swapping algorithm, 
which we will call the ``new'' algorithm.  We compared our results against those of the single-shift version of the 
bulge-chasing algorithm of \cite{KrScWa09}, which we call the ``old'' algorithm.  Our new algorithm can operate directly on the anti-Hessenberg pencil, but the old algorithm requires a preliminary reduction to the partially anti-triangular form
shown in (\ref{eq:partanttri}).  This costs about $n^{2}/8$ moves.  We have not compared against any standard 
non-structure-preserving algorithms, as it was already shown in \cite{KrScWa09} that the structure-preserving 
algorithms are significantly more efficient.  

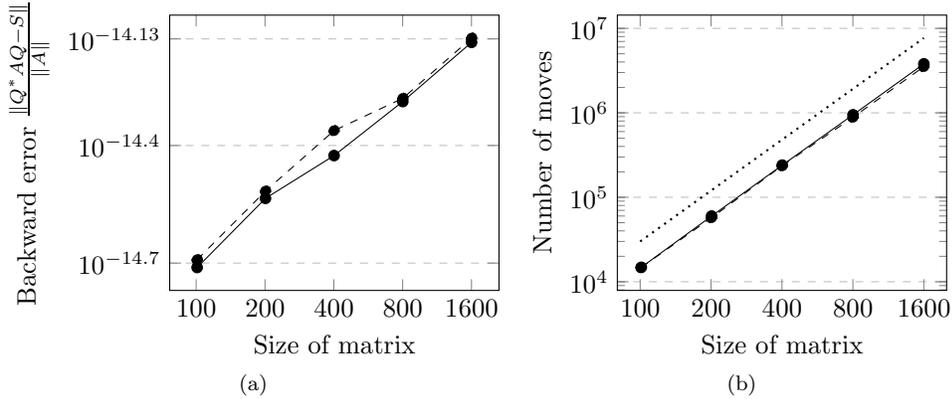
\begin{figure}[ht!]
	\centering
	\subfloat[]{
		\begin{tikzpicture}
		\begin{axis}[
		xlabel={Size of matrix},
		ylabel={Backward error $\frac{\|Q^*AQ - S\|}{\|A\|}$},
		legend pos=north west,
		scale = 0.52,
		ymode=log,
		xmode=log,
		ymajorgrids=true,
		grid style=dashed,
		xtick = {100,200,400,800,1600},
		xticklabels = {100,200,400,800,1600},
		ytick = {2.0E-15, 4.0E-15, 7.5E-15},
		]
		
		\addplot[
		color=black,
		mark=*,
		]
		table [x=n_array, y=backward_errors_pal, col sep=comma] {pal_test_3.csv};
		\addplot[
		color=black,
		mark=*,
		dashed
		]
		table [x=n_array, y=backward_errors_palz, col sep=comma] {pal_test_3.csv};
		
		\end{axis}
		\end{tikzpicture}
	}
	\subfloat[]{
		\begin{tikzpicture}
		\begin{axis}[
		xlabel={Size of matrix},
		ylabel={Number of moves},
		legend pos=north west,
		scale = 0.52,
		ymode=log,
		xmode=log,
		xmin=80,
		xmax=2000,
		ymajorgrids=true,
		grid style=dashed,
		xtick = {100,200,400,800,1600},
		xticklabels = {100,200,400,800,1600},		
		]
		
		\addplot[
		color=black,
		mark=*,
		]
		table [x=n_array, y=swap_count_pal_with_reduction, col sep=comma] {pal_test_3.csv};
		\addplot[
		color=black,
		mark=*,
		dashed
		]
		table [x=n_array, y=swap_count_palz, col sep=comma] {pal_test_3.csv};
		\addplot[
		color=black,
		dotted,
		thick
		]
		coordinates {
			(100, 30000)
			(200, 120000)
			(400, 480000)
			(800, 1920000)
			(1600, 7680000)
		};
		
		\end{axis}
		\end{tikzpicture}
	}
	\caption{Comparison of the old (solid lines) and new (dashed lines) methods for the randomly generated problems.
	For the old algorithm, the move count includes the reduction to partially anti-triangular form.  The dotted line in 
	the move count plot is the function $3n^{2}$, which appears to be parallel to the other two lines. 
	This indicates that the number of moves is of $O(n^2)$, leading to a total flop count of $O(n^3)$.}
	\label{fig:pal-test-3}
\end{figure}

Using the ``old'' and ``new'' methods, we reduced the problem to a pencil $S - \lambda S^{*}$ with $S$ anti-triangular, 
from which we can read the eigenvalues.  Figure \ref{fig:pal-test-3} displays the results.   The left panel shows the
backward errors, which are excellent for both the old and the new methods.  The right panel shows the total number
of moves required as a function of $n$.  We take this as our measure of work.   We see that the two methods 
have very similar performance, with the new method requiring about 5\% fewer moves.  The difference is due 
entirely to the fact that the new method does not require the preliminary reduction.    

The matrices in this test have even dimension, so they use moves of type IIe in the middle.  No refinement 
steps were necessary.   We repeated the experiment with the ten-times stricter tolerance $\epsilon\norm{M}_{F}$
and got nearly identical results.  Refinement steps were occasionally necessary.   

We repeated the experiment with matrices of odd dimension 101, 201, \ldots, 1601, which use moves of type
IIo, and got nearly identical results.  Therefore we have not displayed them.

\begin{figure}[ht!]
	\centering
	\subfloat[]{
		\begin{tikzpicture}
		\begin{axis}[
		xlabel={Size of matrix},
		ylabel={Backward error $\frac{\|Q^*AQ - S\|}{\|A\|}$},
		legend pos=north west,
		scale = 0.52,
		ymode=log,
		xmode=log,
		ymajorgrids=true,
		grid style=dashed,
		xtick = {101,201,401,801,1601},
		xticklabels = {101,201,401,801,1601},
		ytick = {3.2E-15, 6.4E-15, 12.9E-15},
		]
		
		\addplot[
		color=black,
		mark=*,
		]
		table [x=n_array, y=backward_errors_pal, col sep=comma] {pal_test_2.csv};
		\addplot[
		color=black,
		mark=*,
		dashed
		]
		table [x=n_array, y=backward_errors_palz, col sep=comma] {pal_test_2.csv};
		
		\end{axis}
		\end{tikzpicture}
	}
	\subfloat[]{
		\begin{tikzpicture}
		\begin{axis}[
		xlabel={Size of matrix},
		ylabel={Number of moves},
		legend pos=north west,
		scale = 0.52,
		ymode=log,
		xmode=log,
		ymajorgrids=true,
		grid style=dashed,
		xtick = {101,201,401,801,1601},
		xticklabels = {101,201,401,801,1601},		
		]
		
		\addplot[
		color=black,
		mark=*,
		]
		table [x=n_array, y=swap_count_pal_with_reduction, col sep=comma] {pal_test_2.csv};
		\addplot[
		color=black,
		mark=*,
		dashed
		]
		table [x=n_array, y=swap_count_palz, col sep=comma] {pal_test_2.csv};
		\addplot[
		color=black,
		dotted,
		thick
		]
		coordinates {
			(101, 20402)
			(201, 80802)
			(401, 321602)
			(801, 1283202)
			(1601, 5126402)
		};
		
		\end{axis}
		\end{tikzpicture}
	}
	\caption{Comparison of the old (solid lines) and new (dashed lines) methods for the 1D-Laplace boundary control problem \cite[Ex.~2]{KrScWa09}. The dotted line in the number of moves plot is the function $2n^2$.}
	\label{fig:pal-test-2}
\end{figure}
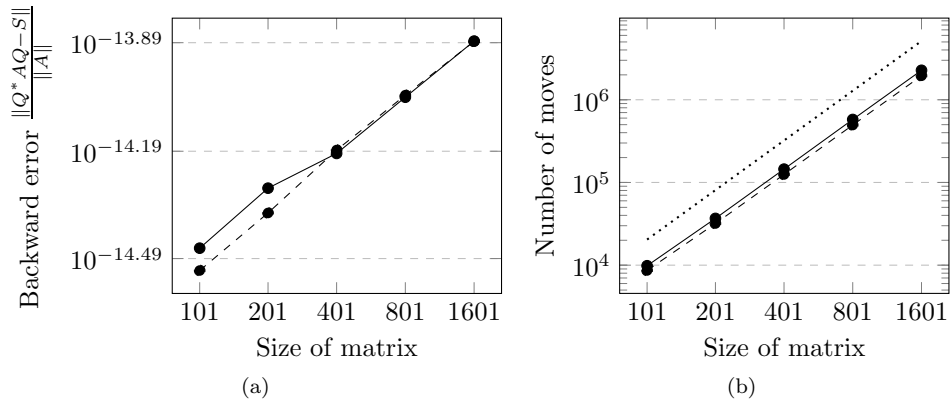

Our second test problem was the discretized 1D-Laplacian system of \cite[Ex.~3]{KrScWa09}.  
In this scalable problem the matrices always have odd dimension, so moves of type IIo are used.  
The results are shown in Figure~\ref{fig:pal-test-2}.   Again we see that both methods have 
similar and excellent backward errors, and the new method requires slightly fewer moves 
than the old method.   This is again due to the fact that the new method does not require a preliminary reduction.  
In this case the difference is about 15\%.    It is larger in this example than it was in the previous one because fewer 
iterations are required in the iterative phase.  In no cases were refinements necessary, even when the stricter tolerance
$\epsilon\norm{M}_{F}$ was enforced.

\subsection*{Tests of moves of types IIo and IIe.}
Since the examples we have considered so far provided barely any test of the refinement procedure, 
we performed stress tests of the moves of types IIo and IIe under extreme conditions.  Type IIo moves 
were tested on matrices of the form 
\begin{equation}
	\begin{bmatrix}
		0 & \ a\\
		a(1+g) & \ c
	\end{bmatrix}
\end{equation}
and type IIe moves were tested on matrices of the form
\begin{equation}
\begin{bmatrix}
0  & \ 0 & \quad a\\
0 & \ b & \quad c\\
a(1+g) & \ d & \quad e
\end{bmatrix}
\end{equation}
where $a$, $b$, $c$, $d$ and $e$ are random numbers generated as $s_110^{t_1} + s_2i10^{t_2}$, with $t_1$ and $t_2$ uniformly distributed in the interval $[-15,0]$ and $s_1$ and $s_2$ random sign bits. This gives numbers with 
wildly varying magnitudes from about $10^{-15}$ up to $10^{0}$.  The parameter $g$ is a positive real number 
that controls the relative gap between the eigenvalues.  We tested both large and small values of $g$, expecting 
that we might sometimes have difficulties when $g$ was very small.  Specifically, we selected four intervals for 
$g$ and generated $10^{5}$ test matrices per interval, with $g$ logarithmically distributed.  We recorded how often 
and how many refinements are required to achieve the tolerance $10\epsilon\norm{M}$. 
The results are shown in Table \ref{table:typeIItests}. 
As expected, refinement steps are never needed when $g$ is large.   Even for very small $g$, the average number 
of refinements is quite small, indicating that refinement steps are very seldom necessary.  Occasionally more than
one refinement step is needed.   There were a few examples with $g$ smaller than $10^{-12}$ where 10 refinements
were required, which is an indication of failure.   In the full algorithm the gap $g$ is controlled by the shifting strategy.
If shifts that are extremely close to the unit circle are avoided, such cases will never occur.  

\begin{table}[]
	\centering
	\begin{tabular}{@{}lllll@{}}
		\toprule
		$g$                    & IIe (average) & IIe (max) & IIo (average) & IIo (max) \\ \midrule
		$[10^{-15}, 10^{-12}]$ & 0.00502    & 10        & 0.08699    & 10        \\
		$[10^{-12}, 10^{-9}]$  & 0.01004    & 3         & 0.089      & 3         \\
		$[10^{-9}, 10^{0}]$    & 0.01413    & 2         & 0.06537    & 2         \\
		$[10^{0}, 10^{15}]$    & 0          & 0         & 0          & 0         \\ \bottomrule
	\end{tabular}
	\caption{The number of refinement steps required to reach the selected tolerance. The columns denoted `average' indicate the average number of refinement steps and the columns denoted `max' indicate the maximum number of refinement steps that was required. Note that because our implementation limits the number of refinement steps to 10, moves where that amount of steps was required are not necessarily accurate.}
	\label{table:typeIItests}
\end{table}

\section{Conclusions}

We have shown that the concept of pole-swapping algorithms, which is a generalization of bulge-chasing algorithms
for the generalized eigenvalue problem, can be extended to palindromic and alternating eigenvalue problems, which
arise in control theory.   We have also introduced a refinement step to make the algorithms (including the algorithms 
in \cite{KrScWa09}) more robust.   Numerical tests of the palindromic case indicate that our pole-swapping algorithm
works well and is slightly faster than the corresponding bulge-chasing algorithm.  


\end{document}